\documentclass[11pt]{article}
\usepackage{epsfig,makeidx,amsmath,amsfonts,latexsym,subfigure}

\newcommand{\E}{{\mathbb E}}

\newcommand {\PP}{\mathbb P}

\begin{document}

\title{Epidemics and vaccination on weighted graphs}

\author{Maria Deijfen \thanks{Department of Mathematics, Stockholm University, 106 91 Stockholm. \newline Email: mia@math.su.se.}}

\date{November 2010}

\maketitle

\thispagestyle{empty}

\begin{abstract}
\noindent A Reed-Frost epidemic with inhomogeneous infection probabilities on a graph with prescribed degree distribution is studied. Each edge $(u,v)$ in the graph is equipped with two weights $W_{(u,v)}$ and $W_{(v,u)}$ that represent the (subjective) strength of the connection and determine the probability that $u$ infects $v$ in case $u$ is infected and vice versa. Expressions for the epidemic threshold are derived for i.i.d.\ weights and for weights that are functions of the degrees. For i.i.d.\ weights, a variation of the so called acquaintance vaccination strategy is analyzed where vertices are chosen randomly and neighbors of these vertices with large edge weights are vaccinated. This strategy is shown to outperform the strategy where the neighbors are chosen randomly in the sense that the basic reproduction number is smaller for a given vaccination coverage.

\noindent
\vspace{0.5cm}

\vspace{0.5cm}

\end{abstract}

\section{Introduction}

The Reed-Frost model is one of the simplest stochastic epidemic models. It was formulated by Lowell Reed and Wade Frost in 1928 (in unpublished work) and describes the evolution of an infection in generations. Each infected individual in generation $t$ ($t=1,2,\ldots$) independently infects each susceptible individual in the population with some probability $p$. The individuals that become infected by the individuals in generation $t$ then constitute generation $t+1$ and the individuals in generation $t$ are removed from the epidemic process. See \cite{vBML} for a description of the asymptotic (as the population size grows to infinity) behavior of the process.

In the original version, an infective individual infects each susceptible individual in the population with the same probability. Realistically however an infective individual has the possibility to infect only those individuals with whom she actually has some kind of social contact. The Reed-Frost model is easily modified to capture this by introducing a graph to represent the social structure in the population and then let the infection spread on this graph. More precisely, an infective individual infects each neighbor in the graph independently with some probability $p$.

When analyzing epidemics on graphs, the graph is usually taken to be unweighted with respect to the infection, that is, transmission takes place along all edges with the same probability. In this paper however, inhomogeneity will be incorporated in the transmission probability by aid of weights on the edges. More precisely, each edge $(u,v)$ in the graph is assigned two weights $W_{(u,v)}$ and $W_{(v,u)}$ that are assumed to take values in [0,1]. The probability that $u$ infects $v$ if $u$ gets infected is then given by $W_{(u,v)}$ and vice versa. Note that it may well be that $W_{(u,v)}\neq W_{(v,u)}$. We shall mainly consider i.i.d.\ weights, although we briefly treat weights that are determined by the degrees of the vertices in Section 3.

To describe the underlying network, we shall use the so called configuration model \cite{MR-95,MR-98}. Once the graph has been generated, each edge is equipped with two weights as described above. Basically, the configuration model takes a probability distribution with support on positive integers as input and generates a graph with this particular degree distribution; see Section 2 for further details. The degree distribution is indeed an important characteristic of a network with a large impact on the properties of the network and it is therefore desirable to be able to control this in a graph model. Furthermore, the configuration model exhibits short distances between the vertices, which is in agreement with empirical findings; see \cite{vHHvM}. Epidemics on un-weighted graphs generated by the configuration model has previously been studied in \cite{A-98,A-99,BJML}. Related results have also appeared in the physics literature \cite{N}.

An important quantity in epidemic modeling is the epidemic threshold, commonly denoted by $R_0$. It is defined as a function of the parameters of the model such that a large outbreak in the epidemic has positive probability if and only if $R_0>1$. Expressions for $R_0$ typically stems from branching process approximations of the initial stages of the epidemic. It is well-known from branching process theory that the process has a positive probability of exploding if and only if the expected number of children of an individual exceeds 1. A natural candidate for $R_0$ is hence the expected number of new cases caused by a typical infective in the beginning of the time course. For this reason, the epidemic threshold is often referred to as the basic reproduction number.

The main goal of the paper is to study how the epidemic threshold is affected by vaccination strategies based on the edges weights. To this end, we assume that a perfect vaccine is available that completely removes vaccinated individuals from the epidemic process. The simplest possible vaccination strategy, usually referred to as random vaccination, is to draw a random sample from the population and then vaccinate the corresponding individuals. An alternative, known as acquaintance vaccination, is to choose individuals randomly and then, for each chosen individual, vaccinate a random neighbor rather than the individual itself \cite{CHbA,BJML}. The idea is that, by doing this, individuals with larger degrees are vaccinated. We shall study a version of acquaintance vaccination where, instead of vaccinating a random neighbor, the neighbor with the largest weight on its edge from the sampled vertex is vaccinated. In a human population, this correspond to asking individuals to name their \emph{closest} friend (in some respect) instead of just naming a random friend. It is demonstrated that this is more efficient than standard acquaintance vaccination, in the sense that the basic reproduction number with the weight based strategy is smaller for a given vaccination coverage.

Throughout this paper we shall use the term ``infection'' to refer to the phenomenon that is spreading on the network. We remark that this does not necessarily consist of an infectious disease spreading in a human population, but may also refer to other infectious phenomena such as a computer virus spreading in a computer network, information routed in a communication net or a rumor growing in a social media. In many of these situations the connections are indeed highly inhomogeneous. Furthermore, depending on what type of spreading phenomenon that is at hand, the term vaccination can refer to different types of immunization.

Epidemics on weighted graphs have been very little studied so far and there are few theoretical results. See however \cite{N} for an approach based on generating function and \cite{Gang,Schumm} for simulation studies. We mention also the recent work on first passage percolation on random graphs by Bhamidi et al.\ \cite{BvHH:1,BvHH:2,BvHH:3}. There, each edge in a graph generated according to the configuration model is equipped with an exponential weight and the length and weight of the weight-minimizing path between two vertices are studied. Interpreting the weights as the traversal times for an infection, this can be related to the time-dynamics of an epidemic.

The rest of the paper is organized so that the graph model and the epidemic model are described in more detail in Section 2. In Section 3, expressions for the epidemic thresholds are given and calculated for some specific weight distributions. Section 4 is devoted to vaccination: In Section 4.1, a weight based acquaintance vaccination strategy for weights with a continuous distribution is described and an expression for the epidemic threshold is derived. Section 4.2 treats a strategy for a two-point weight distribution. The findings are summarized in Section 5, where also some directions for further work are given. We shall throughout refrain from giving rigorous details for the underlying branching process approximations, but instead focus on heuristic derivations of the epidemic quantities. Indeed, what needs to be proved is basically that the branching process approximations hold long enough so that conclusions for the branching processes are valid also for the epidemic processes. This however is not affected by weights on edges (as long as these are not functions of the structure of the graph) and hence rigorous details can presumably be filled in by straightforward modifications of the arguments in \cite{BJML} (the degree based weights mentioned in Section 3 might however require some more work).

\section{Description of the model}

We consider a population of size $n$ represented by $n$ vertices. The graph representing the connections in the population is generated by the configuration model. To produce the graph, a probability distribution with support on the non-negative integers is fixed and each vertex $u$ is independently equipped with a random number of half-edges $D_u$ according to this distribution. These half-edges are then paired randomly to create the edges in the graph, that is, first two half-edges are picked at random and joined, then another two half-edges are picked at random from the set of remaining half-edges and joined, etc. If the total number of half-edges is odd, a half-edge is added at a randomly chosen vertex to pair with the last half-edge.

This procedure gives a multi-graph, that is, a graph where self-loops and multiple edges between vertices may occur. If $D$ has finite second moment however, there will not be very many of these imperfections. In particular, the probability that the resulting graph is simple will be bounded away from 0 as $n\to\infty$; see \cite[Lemma 5.5]{BJML} or \cite[Theorem 7.10]{vH}. If $D$ has finite second moment we can hence condition on the event that the graph is simple, and work under this assumption. Another option is to erase self-loops and merge multiple edges, which asymptotically does not affect the degree distribution if $D$ has finite second moment; see \cite[Theorem 7.9]{vH}. Henceforth we shall hence assume that $D$ has finite second moment and ignore self-loops and multiple edges.

When the graph has been generated, each edge $(u,v)$ is assigned two weights $W_{(u,v)}$ and $W_{(v,u)}$ that are assumed to take values in [0,1]. This can be thought of as if each one of the half-edges that is used to create the edge independently receives a weight. The epidemic spread is initiated in that one randomly chosen vertex is infected. This vertex constitutes generation 1. The epidemic then propagates in that each vertex $u$ in generation $t$ ($t=1,2,\ldots$) infects each susceptible neighbor $v$ independently with probability $W_{(u,v)}$. Generation $t+1$ then consists of the vertices that are infected by the vertices in generation $t$ and the vertices in generation $t$ are removed from the epidemic process.

We shall mainly restrict to the case where the weights are taken to be independent. However, we mention also the possibility to let them be functions of the degrees of the vertices:\medskip

\noindent\textbf{Independent weights.} The weights are taken to be i.i.d.\ copies of a random variable $W$ that takes values in [0,1]. The distribution of $W$ can be defined in many different ways:

\begin{itemize}
\item[$\bullet$] As an intrinsic distribution on [0,1], for instance a uniform distribution or, more generally, a Beta distribution.
\item[$\bullet$] By letting $N$ be an integer valued random variable, indicating for instance how many times a given vertex contacts a given neighbor during some time period, and then setting $W\stackrel{d}{=}1-(1-p)^N$, with $p\in[0,1]$ denoting the probability of infection at a given contact.
\item[$\bullet$] By, similarly, letting $X$ be a positive random variable, interpreted as the (subjective) strength of a connection, and then for instance setting $W\stackrel{d}{=}\mathbf{1}\{X\geq \theta\}$ for some $\theta\geq 0$ or $W\stackrel{d}{=}1-\alpha^X$ for $\alpha\in[0,1]$. Alternatively, $X$ could be interpreted as the resistance involved in a connection and $W$ modeled as a decreasing function of $X$.
\end{itemize}

\noindent \textbf{Degree dependent weights.} The weights of an edge $(u,v)$ could also be modeled as functions of $D_u$ and $D_v$. We shall consider the case when $W_{(u,v)}=g(D_u)$ for some function $g$ that takes values in [0,1]. All outgoing edges from $u$ hence have the same weight, and independent trials with this success probability determine whether the edges are used to transmit infection. With $g$ increasing, this setup means that vertices with large degree have a larger probability of infecting their neighbors, for instance in that they tend to be more active. With $g$ decreasing, high degree vertices are instead less likely to infect their neighbors, which might be the case for instance in a situation where high degree vertices have weaker bonds to their acquaintances.

\section{Epidemic threshold}

As mentioned in the introduction, expressions for epidemic thresholds usually come from branching process approximations of the initial stages of an epidemic. As for epidemics on graphs, branching process approximations are typically in force as soon as the graph is tree-like, that is, if with high probability the graph does not contain short cycles. This means that the neighbors of a given infective in the beginning of the time course are susceptible with high probability and hence the initial stages of the generation process of infectives is well approximated by a branching process. Under the assumption that the degree distribution has finite second moment, the configuration model is indeed tree-like, allowing for such an approximation; see e.g.\ \cite{vHHvM,BJML} for details. The epidemic threshold is then given by the reproduction mean in the approximating branching process, which in turn is given by the expected number of new cases generated by an infective vertex in the beginning of the epidemic. When calculating this, one should not consider the initial infective, since this vertex might be atypical, but rather an infective vertex in, say, the second generation.

Let $\{p_k\}_{k\geq 0}$ be the probabilities defining the degree distribution in the configuration model. Then the initial infective has degree distribution $\{p_k\}$, while the neighbors of this vertex have the size biased degree distribution $\{\tilde{p}_k\}$ defined by
$$
\tilde{p}_k=\frac{kp_k}{\mu},
$$
where $\mu=\sum kp_k$ denotes the mean degree. The infective vertices in the second (and later) generations hence have degree distribution $\{\tilde{p}_k\}$. Denote by $\widetilde{D}$ a random variable with this distribution.\medskip

\noindent\textbf{Independent weights.} Consider an infected vertex in the second generation. One neighbor of this vertex must have transmitted the infection and can hence not get reinfected, while the other neighbors are with high probability susceptible. The number of new cases generated by the vertex is hence distributed as
\begin{equation}\label{eq:new_cases}
\sum_{i=1}^{\tilde{D}-1}\mathbf{1}\{\mbox{neighbor $i$ infected}\}.
\end{equation}
If the weights are i.i.d.\ copies of $W$, then the mean of the indicators is $\gamma:=\E[W]$ and we get
$$
R_0=\gamma\E[\widetilde{D}-1]=\gamma\left(\mu+\frac{\mbox{Var}(D)-\mu}{\mu}\right).
$$
In this case the epidemic threshold is hence the same as in a model with constant infection probability $p=\gamma$ ; see \cite{A-99,BJML}. Note that, for degree distributions with large variance, $R_0$ can be large even if $\mu$ is small. Also note that the above reasoning remains valid in a situation where the weights $W_{(u,v)}$ and $W_{(v,u)}$ on a given edge are correlated, as long as the weights are independent between edges. In fact, as long as the transmission between separate links are i.i.d., the whole epidemic process is equivalent to a Reed Frost model (on the configuration model) with $p=\E[W]$.\medskip

\noindent \textbf{Degree dependent weights.} Assume that $W_{(u,v)}=g(D_u)$. Since $W_{(u,v)}$ does not depend on $D_v$, the degree distribution of an infective in the second generation is $\{\tilde{p}_k\}$. Conditionally on its degree $\widetilde{D}=\tilde{d}$, the number of new cases generated by an infective in the second generation is Bin($\tilde{d}-1,g(\tilde{d})$)-distributed. It follows that
$$
R_0^{deg}=\E[(\widetilde{D}-1)g(\widetilde{D})].
$$
Let $R_0^{h1}$ denote the basic reproduction number for an epidemic with a homogeneous infection probability given by the transmission probability $\E[g(\tilde{D})]$ for a randomly chosen half-edge, that is,
$$
R^{h1}_0=\E[g(\widetilde{D})]\E[\widetilde{D}-1].
$$
When $g$ is an increasing function we have $R_0^{deg}\geq R_0^{h1}$, due to the positive correlation between $g(\widetilde{D})$ and $\widetilde{D}$, while, if $g$ is decreasing, then $R_0^{deg}\leq R_0^{h1}$. Another comparison that might be relevant is to relate $R_0^{deg}$ to the basic reproduction number for an epidemic with a homogeneous infection probability given by $\E[g(D)]$, that is, an epidemic where the infection probability $g(D)$ for a vertex with degree $D$ is averaged over all possible degrees. The basic reproduction number in such an epidemic is given by
$$
R_0^{h2}=\E[\widetilde{D}-1]\E[g(D)].
$$

\noindent \textbf{Example 3.1.} First take $D\sim$ Po($\mu$). It is not hard to see that then $\tilde{p}_k=p_{k-1}$. Take $g(x)=\mathbf{1}\{x\geq \theta\}$. If $W_{(u,v)}=g(D_u)$, this means that only vertices with degree at least $\theta$ transmit the infection. We have
$$
\E[\mathbf{1}\{\widetilde{D}\geq \theta\}]=\PP(D\geq \theta-1)
$$
and
$$
\E[(\widetilde{D}-1)\mathbf{1}\{\widetilde{D}\geq \theta\}]=\sum_{k\geq \theta}(k-1)\frac{kp_k}{\mu}=\sum_{k\geq \theta}\frac{\mu^{k-1}}{(k-2)!}e^{-\mu}=\mu\PP(D\geq \theta-2).
$$
Hence
$$
\begin{array}{lcl}
R_0^{h2} & = & \mu\PP(D\geq \theta)\\
R_0^{h1} & = & \mu\PP(D\geq \theta-1)\\
R_0^{deg} & = & \mu\PP(D\geq \theta-2).
\end{array}
$$
With $g(x)=\alpha^x$ for $\alpha\in(0,1)$, we get
$$
\E[\alpha^D]=e^{-\mu(1-\alpha)}\quad\mbox{and}\quad\E[\alpha^{\widetilde{D}}]=\alpha e^{-\mu(1-\alpha)}.
$$
Furthermore
$$
\E[(\widetilde{D}-1)\alpha^{\widetilde{D}}]=\sum_{k\geq 1}(k-1)\alpha^k\frac{kp_k}{\mu}=
\alpha^2\sum_{k\geq 0}\alpha^kp_k=\alpha^2e^{-\mu(1-\alpha)}.
$$
Hence
$$
\begin{array}{lcl}
R_0^{h2} & = & \mu e^{-\mu(1-\alpha)}\\
R_0^{h1} & = & \mu\alpha e^{-\mu(1-\alpha)}\\
R_0^{deg} & = & \mu\alpha^2e^{-\mu(1-\alpha)}.
\end{array}
$$
\hfill$\Box$\medskip

\noindent \textbf{Example 3.2.} Now take a distribution with $p_k\sim ck^{-3.5}$. In this case exact computations are out of reach but numerical values of the thresholds are easily obtained. We give an example with $g(x)= x^{-\tau}$ for $\tau\in[0,1]$. The initial degrees in the graph have been modified to give a mean of 4. In Figure 1, the basic reproduction numbers are plotted against $\tau$, showing that $R_0^{h2}>R_0^{h1}>R_0^{deg}$. The homogeneous epidemics remain supercritical (that is, their reproduction numbers exceed 1) at $\tau=1$, while the epidemic with degree dependent weights becomes subcritical for $\tau$ close to 1. Indeed, the degree dependent epidemic is subcritical at $\tau=1$ for any degree distribution.\hfill$\Box$\medskip

\begin{figure}
\begin{center}
\mbox{\epsfig{file=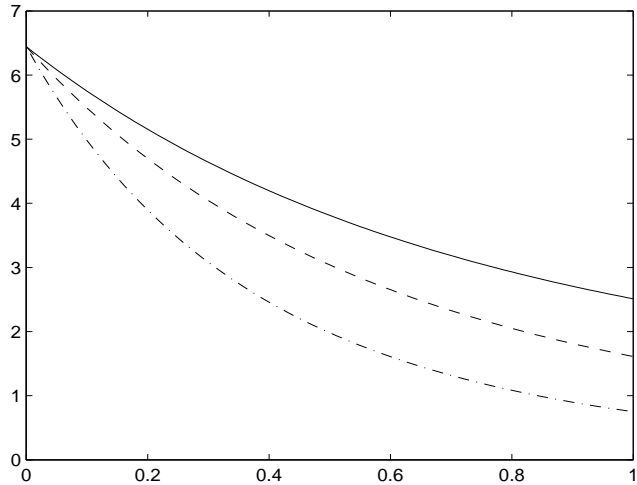,width=0.8\textwidth,
height=0.35\textheight}}
\end{center}
\caption{The basic reproduction numbers $R_0^{h2}$ (solid line), $R_0^{h1}$ (dashed line) and $R_0^{deg}$ (dash-dotted line) with $g(x)=x^{-\tau}$ plotted against $\tau$. The degree distribution is a power law with exponent 3.5 and mean 4.}
\end{figure}

For the remainder of the paper we shall restrict to the case with independent weights.

\section{Vaccination}

We now proceed to analyze a version of the so called acquaintance vaccination strategy. To this end, suppose that a perfect vaccine is available that prevents vaccinated vertices from participating in the epidemic process and that this vaccine is distributed prior to the start of the epidemic. More precisely, first we generate the underlying graph and assign the edge weights, then we choose which vertices that are to be vaccinated and finally, when the vaccine has been distributed, we analyze the epidemic spread among the unvaccinated vertices.

The simplest vaccination scheme is to vaccinate each vertex independently with some probability $v$. We shall refer to this as uniform vaccination and write $R_v^{\verb"U"}$ for the corresponding basic reproduction number. By reasoning as in the case without vaccination and keeping in mind that only unvaccinated vertices can be infected, it is not hard to see that $R_v^{\verb"U"}$  is obtained by multiplying the expression for the case without vaccination with $(1-v)$, that is,
\begin{equation}\label{eq:RvU}
R_v^{\verb"U"}=(1-v)\gamma\left(\mu+\frac{\mbox{Var}(D)-\mu}{\mu}\right).
\end{equation}
An alternative strategy, referred to as acquaintance vaccination, is to vaccinate neighbors of the chosen vertices rather than the vertices themselves; see \cite{CHbA}. More precisely, each vertex is sampled independently with probability $s$ and, for each sampled vertex, a randomly chosen neighbor is vaccinated. A neighbor that is chosen via more than one vertex is (of course) vaccinated only once. The fact that two vertices can both pick a common neighbor to receive vaccination means that the asymptotic fraction of vaccinated vertices $v(s)$, referred to as the \emph{vaccination coverage}, is smaller than $s$. In fact, in many cases it is not possible to push an epidemic below criticality (that is, to obtain a basic reproduction number smaller than 1) even when $s=1$.  This motivates a strategy where a vertex can be sampled more than once and thereby have more than one neighbor vaccinated.

In \cite{BJML}, a strategy is analyzed where each vertex is sampled independently a Po($\beta$) number of times and each time a vertex is sampled, a randomly chosen neighbor is vaccinated. When there are weights on the edges however, more efficient strategies are possible. Here we shall consider strategies where, rather than choosing neighbors randomly for vaccination, neighbors with large weights on their edges from the sampled vertices are chosen. We treat the case with independent directed weights. In Section 4.1, the weights are assumed to come from a continuous distribution and Section 4.2 is devoted to a strategy for two-point distributions. 

\subsection{Weight based acquaintance vaccination: continuous weights}

Assume that the directed edge weights are i.i.d.\ realizations from a continuous probability distribution on [0,1]. Each vertex $u$ is sampled independently a Po($\beta$) number of times. Write $\mathcal{N}_u$ for the set of neighbors of a vertex $u$ and, for $i=1,\ldots,|\mathcal{N}_u|$, let $v_i\in\mathcal{N}_u$ be the vertex corresponding to the $i$:th largest element in $\{W_{(u,v)}:v\in\mathcal{N}_u\}$ (note that $v_i$ is almost surely unique, since the weights are assumed to come from a continuous distribution). Then, if $u$ is sampled $i\leq D_u$ times, the neighbors $v_1,\ldots,v_i$ are vaccinated. If $i\geq D_u$, all neighbors of $u$ are vaccinated and, if $D_u=0$, no action is taken. This will be referred to as weight based acquaintance vaccination.

To derive an expression for the vaccinated fraction $v(\beta)$ of the population, let $V\sim$ Po$(\beta)$ represent the number of times that a given vertex is sampled. The probability that a randomly chosen vertex $u$ is not chosen for vaccination by a given neighbor with degree $k$ is given by
\begin{equation}\label{eq:r_k}
r_k=\sum_{i=0}^{k-1}\PP(V=i)\left(1-\frac{i}{k}\right).
\end{equation}
Since the neighbors of $u$ have degree distribution $\{\tilde{p}_k\}$, the probability that $u$ avoids being chosen for vaccination by a given neighbor equals
$$
\alpha=\sum_{k\geq 1}r_k\tilde{p}_k.
$$
If $u$ has degree $j$, then the probability that $u$ is not vaccinated is $\alpha^j$. The degree distribution of $u$ is $\{p_j\}$ and, since the fraction of unvaccinated vertices coincides with the probability that a randomly chosen vertex is not vaccinated, we obtain $v(\beta)$ from the equation
$$
1-v(\beta)=\sum_{j\geq 0}\alpha^jp_j.
$$
Next, to identify the epidemic threshold in a population vaccinated according to the weight based acquaintance strategy, we shall employ a branching process approximation of the initial stages of the epidemic. The process however is slightly more complicated than in the case with uniform vaccination, and this is because the knowledge that an edge has not been used for vaccination carries information of the degrees of the corresponding vertices. The process is analogous to the one used in \cite{BJML}. To describe it, say that a directed edge $(u,w)$ is \emph{used for vaccination} if $u$ is sampled and chooses the neighbor $w$ for vaccination. Furthermore, conditionally on the weight, a directed edge is said to be \emph{open for transmission} if a bernoulli trial with success probability given by the weight of the edge results in a success. A directed edge that is not used for vaccination and that is open for transmission is called \emph{dangerous}.

An ``individual'' in the branching process now consists of an unvaccinated vertex $u$ along with a dangerous outgoing edge $(u,w)$. The individual then gives birth to a new individual if the vertex $w$ is unvaccinated and in turn has a dangerous edge $(w,w')$ pointing out from it. Note that an unvaccinated vertex can hence give rise to several individuals (if it has more than one outgoing dangerous edge) or no individuals at all (if it does not have any outgoing dangerous edges). Furthermore, the individuals reproduce independently.

It is not hard to see that the epidemic has a positive probability of taking off if and only if the above branching process has a positive probability of exploding: With positive probability the initial infective in unvaccinated and with positive probability it has at least one dangerous out-edge. The propagation of the epidemic from the vertices that are hit by these dangerous edges is then approximated by the above branching process. To find an expression for the reproduction mean of the process, consider a given unvaccinated vertex $u$ along with an outgoing dangerous edge $(u,w)$. How many new individuals does this give rise to? First, the degree distribution of $w$, which is size biased, is now affected also by the information that the edge $(w,u)$ has not been used for vaccination. Write $A$ for the latter event and $\PP_A$ and $\E_A$ for probability and expectation respectively conditionally on $A$. We get
$$
\PP_A(D_w=k)=\frac{r_k\tilde{p}_k}{\sum r_k\tilde{p}_k}=\frac{r_k\tilde{p}_k}{\alpha},
$$
where $r_k$ is defined in (\ref{eq:r_k}). If $w$ has degree $k$, then the probability that $w$ is not chosen for vaccination by any of its other $k-1$ neighbors (apart from $u$ that, by assumption, does not have a dangerous edge to $w$) is $\alpha^{k-1}$. Write $H$ for the number of dangerous edges $(w,w')$ with $w'\neq u$. Note that, conditionally on the degree of $w$, the event that $w$ is unvaccinated (which carries information on the in-weights of $w$) does not affect $H$ (which is determined by the out-weights of $w$). The reproduction mean of the branching process is hence given by
\begin{equation}\label{eq:RbW}
R_\beta^{\verb"W"}=\sum_{k\geq 2}\PP_A(D_w=k)\alpha^{k-1}\E_{A,k}[H],
\end{equation}
where $\E_{A,k}[H]:=\E_A[H|D_w=k]$ (below we use $\PP_{A,k}$ to denote the corresponding probability and $\PP_k$ to denote probability conditional only on that $D_w=k$). It remains to quantify this expectation. Clearly $H$ is affected by the number of times $V_w$ that $w$ is sampled to name a neighbor in the vaccination procedure, which in turn is affected by the information that $(w,u)$ was not used for vaccination. Specifically, when $D_w=k,$ we have, for $i=0,\ldots,k-1$, that
\begin{eqnarray*}
\PP_{A,k}(V_w=i) & = & \frac{\PP_k(A|V_w=i)\PP_k(V_w=i)}{\PP_k(A)}\\
& =& \frac{\left(1-\frac{i}{k}\right)\PP(V_w=i)}{r_k}.
\end{eqnarray*}
Let $W^{(k)}_j$ denote a random variable distributed as the $j$:th smallest in a collection of $k$ independent weight variables. If $V_w=i$ ($i=1,\ldots,k-1$), then the out-edges with the $i$ largest weights are used for vaccination. The remaining $k-i$ out-edges are dangerous with a probability given by the expectation of their weights. Note however that we do not want to count the edge to $u$, whose weight indeed belongs to the $k-i$ smallest since, by assumption, it is not used for vaccination. The ordering of the weight on the edge to $u$ among the remaining $k-i$ out-weights is uniform on $\{1,\ldots,k-i\}$. We obtain
\begin{equation}\label{eq:vvHk}
\E_{A,k}[H|V_w=i]=\left(1-\frac{1}{k-i}\right)\sum_{j=1}^{k-i}\E[W^{(k)}_j].
\end{equation}
Note that, if $V_w=0$, then each one of the $k-1$ out-edges from $w$ to $\mathcal{N}_w\backslash \{u\}$ is dangerous independently with probability $\gamma$, and the above expression reduces to $(k-1)\gamma$. If $V_w\geq k-1$, then all out-edges (except $(w,u)$) are used for vaccination meaning that there are no dangerous out-edges. Hence
$$
\E_{A,k}[H]=\sum_{i=0}^{k-2}\PP_{A,k}(V_w=i)\E_{A,k}[H|V_w=i].
$$
This concludes the derivation of the reproduction mean (\ref{eq:RbW}).

Calculating the reproduction mean involves calculating expectations of order statistics (c.f.\ (\ref{eq:vvHk})). Finding analytical expressions for such expectations is typically not possible. However, the density of the $j$:th smallest observation in a collection of $k$ i.i.d.\ variables with density $f$ and distribution function $F$ is given by
\begin{equation}\label{eq:ord_stat}
f_{k,j}(x)=\frac{\Gamma(k+1)}{\Gamma(j)\Gamma(k+1-j)}(F(x))^{j-1}(1-F(x))^{k-j}f(x),
\end{equation}
where $\Gamma(\cdot)$ denotes the gamma function. For a given weight distribution $F$, the mean can hence be calculated by aid of numerical integration. A particularly easy case is when the weights are uniform on [0,1]. Then
$$
W_j^{(k)}\sim\mbox{Beta}(j,k+1-j)
$$
so that $\E[W_j^{(k)}]=j/(k+1)$, and hence
\begin{equation}\label{eq:uni_dan}
\E_{A,k}[H|V_w=i] = \left(1-\frac{1}{k-i}\right)\sum_{j=1}^{k-i}\frac{j}{k+1} = \frac{(k-i-1)(k-i+1)}{2(k+1)}.
\end{equation}

We now want to compare the epidemic threshold for the weight based strategy to the threshold for the standard acquaintance vaccination strategy, where neighbors are chosen randomly. In \cite{BJML}, the vaccination coverage $v(\beta)$ for the latter strategy is shown to be given by
$$
1-v(\beta)=\sum_{j\geq 0}\alpha^jp_j
$$
with $\alpha=\sum e^{-\beta/k}\tilde{p}_k$. Furthermore, for a homogeneous infection probability $p$, the basic reproduction number is shown to be
\begin{equation}\label{eq:RbA}
R_\beta^{\verb"A"}=p\sum_{k\geq 2}(k-1)\alpha^{k-2}e^{-2\beta/k}\tilde{p}_k.
\end{equation}
Straightforward modifications of the arguments leading up to these expressions reveals that they apply also for the inhomogeneous case with independent weights, with $p$ replaced by the mean weight $\gamma$.\medskip

\noindent \textbf{Example 4.1.1.} Let the edge weights be uniformly distributed on $[0,1]$. Then $\gamma=1/2$ and, using (\ref{eq:uni_dan}), the reproduction mean $R_\beta^{\verb"W"}$ in (\ref{eq:RbW}) is easily calculated for a given degree distribution $\{p_k\}$. Figure 2 shows the basic reproduction number $R_\beta^{\verb"W"}$ plotted against the vaccination coverage $v(\beta)$ when the degree distribution is Po(6). The plot also shows the reproduction number for standard acquaintance vaccination and for uniform vaccination. For a given vaccination coverage, we have $R^{\verb"U"}>R^{\verb"A"}>R^{\verb"W"}$, although the difference between standard acquaintance vaccination and the weight based strategy is quite small. Note however that, in practical situations also a small gain could be valuable: The vaccination coverage required to push the epidemic threshold below 1 -- thereby preventing large outbreaks -- is referred to as the \emph{critical vaccination coverage}. Clearly, when fighting an infectious disease in a large human population for instance, even a very small decrease in the critical vaccination coverage might imply large savings in terms of vaccination costs.\hfill$\Box$\medskip

\noindent \textbf{Example 4.1.2.} Let the weights have a Beta distribution with parameters 0.5 and 2.5; see Figure 3. In this case it is not possible to write down analytical expressions in closed form for $R^{\verb"W"} $ but it is easily computed numerically. Figure 4 shows the basic reproduction numbers plotted against the vaccination coverage when the degree distribution is Po(14). In this case the weight based strategy performs clearly better than the standard acquaintance vaccination. In particular, the critical vaccination coverage for uniform vaccination and standard acquaintance vaccination is 0.58 and 0.53 respectively, while for the weight based strategy it is decreased to 0.47. The reason is that the weight distribution is right-skewed: Most weights are small but there is a thick right-tail with large weights, and by getting rid of these large weights the mean in the weight distribution is decreased more than in the uniform case.\hfill$\Box$\medskip

\noindent \textbf{Example 4.1.3.} Finally, let the weights have the same Beta distribution as in the previous example, but take the degree distribution to be a power-law with exponent 3.5 and the same mean 14 as in the Poisson distribution. Figure 5 displays the basic reproduction numbers in this case. Again the weight based strategy performs better than standard acquaintance vaccination. In this case however, the most striking feature is the difference between the uniform vaccination and the acquaintance based strategies: when the degree distribution is a power law, the basic reproduction number is pushed down very effectively by targeting high degree vertices.\hfill$\Box$\medskip

\begin{figure}[p]
\begin{center}
\mbox{\epsfig{file=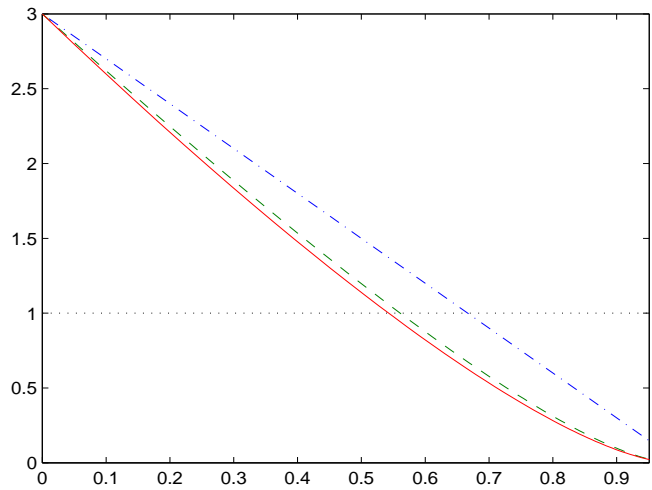,width=0.8\textwidth,
height=0.35\textheight}}
\end{center}
\caption{Basic reproduction numbers with $U(0,1)$ weights for a Po(6) degree distribution plotted against the vaccination coverage: the weight based acquaintance strategy (solid line), the standard acquaintance vaccination (dashed line) and uniform vaccination (dash-dotted line).}
\end{figure}

\begin{figure}[p]
\begin{center}
\mbox{\epsfig{file=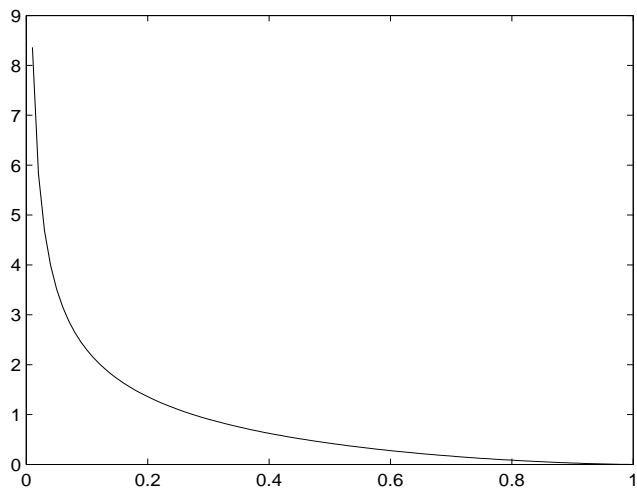,width=0.8\textwidth,
height=0.35\textheight}}
\end{center}
\caption{A Beta density with parameter 0.5 and 2.5.}
\end{figure}

\begin{figure}[p!]
\begin{center}
\mbox{\epsfig{file=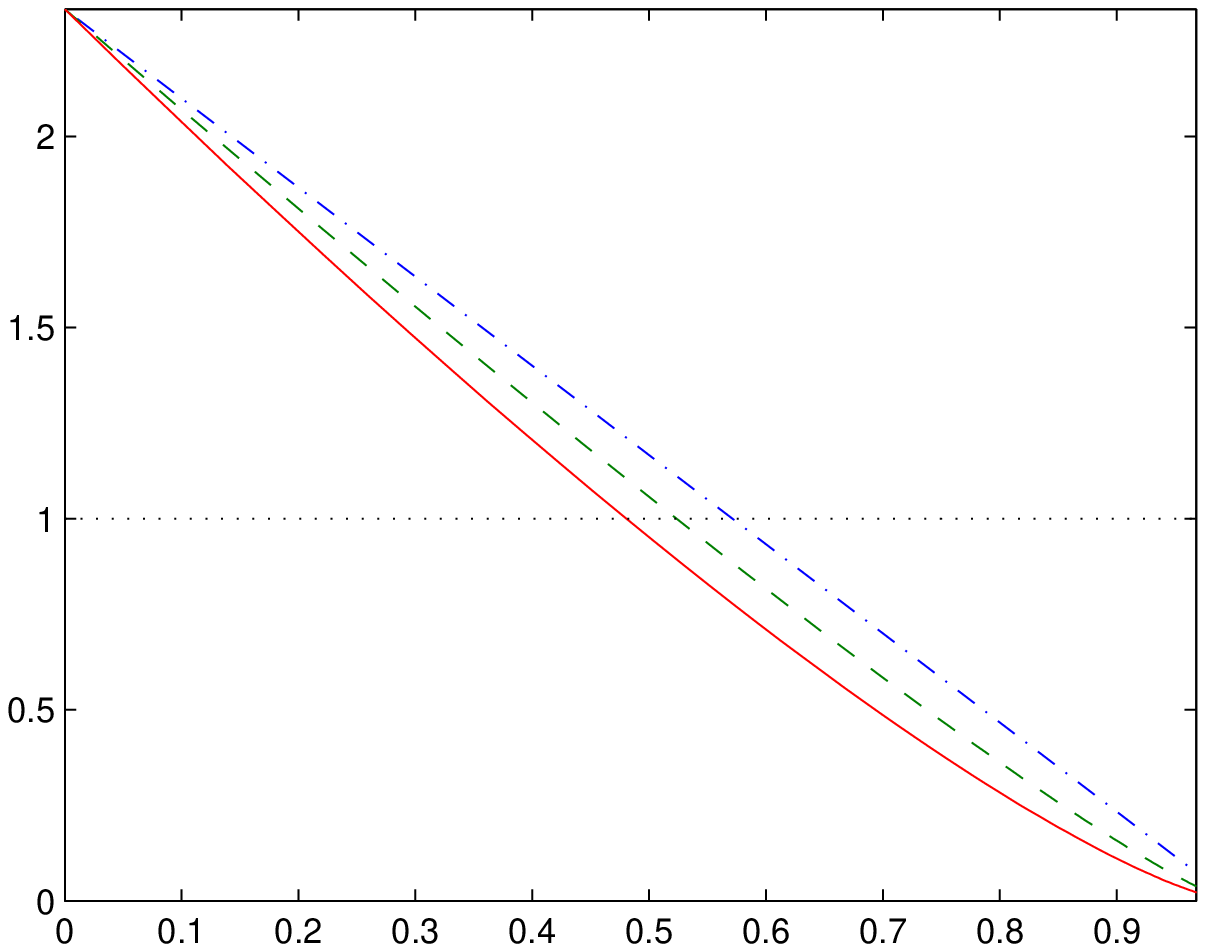,width=0.8\textwidth,
height=0.34\textheight}}
\end{center}
\caption{Basic reproduction numbers with Beta(0.5,2.5) weights for a Po(14) degree distribution plotted against the vaccination coverage: the weight based acquaintance strategy (solid line), the standard acquaintance vaccination (dashed line) and uniform vaccination (dash-dotted line).}
\end{figure}

\begin{figure}[p!]
\begin{center}
\mbox{\epsfig{file=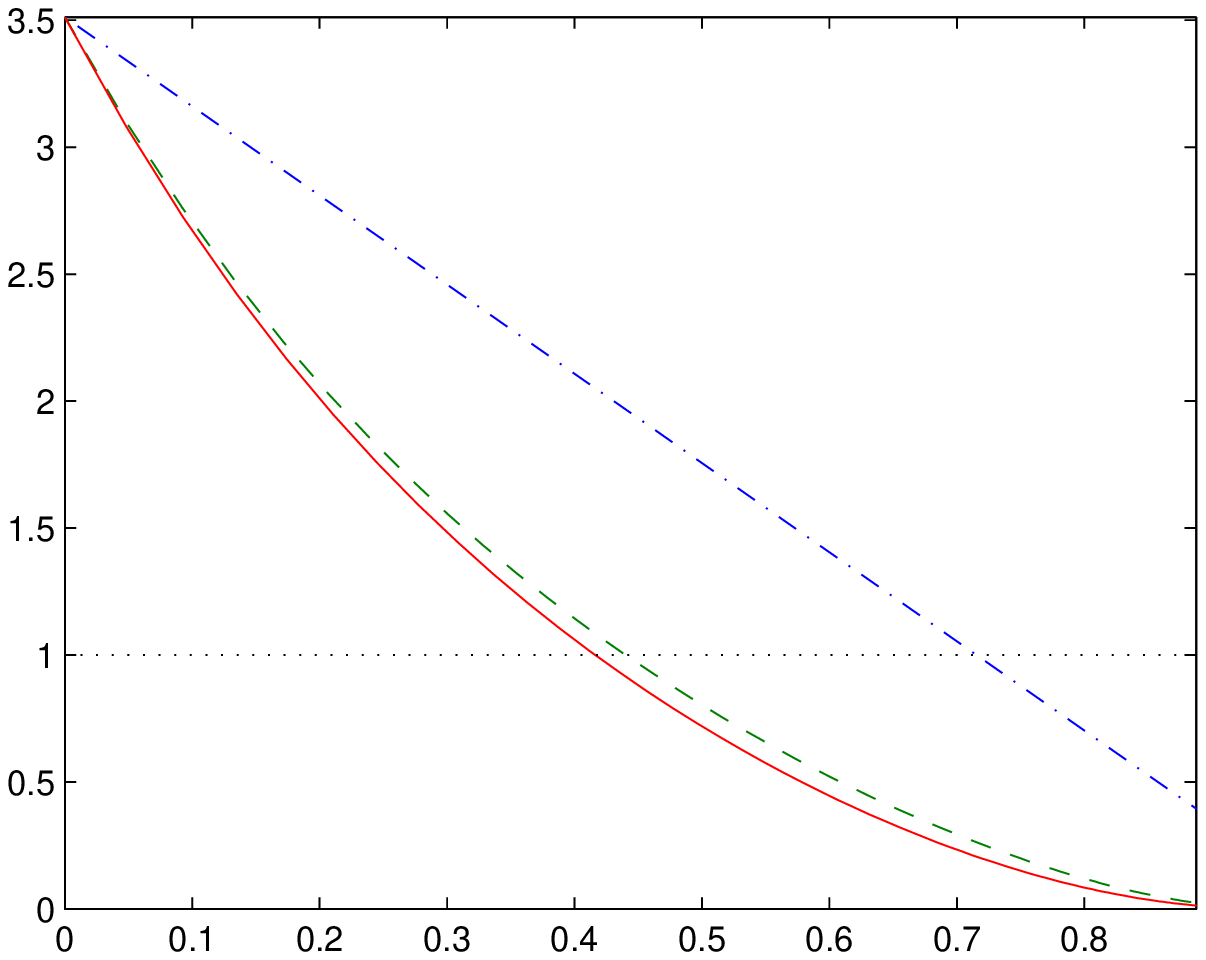,width=0.8\textwidth,
height=0.34\textheight}}
\end{center}
\caption{Basic reproduction numbers with Beta(0.5,2.5) weights for power-law degree distribution with exponent 3.5 and mean 14 plotted against the vaccination coverage.}
\end{figure}

\subsection{Weight based vaccination: two-point weights}

The finding in the previous section that the weight based strategy performs well for right-skewed weight distributions in the continuous case might lead one to suspect that the strategy is particularly useful for a ``polarized'' discrete distribution. In this section we analyze the simple case when the weights have a two-point distribution. As a motivation we can think of a network having two types of directed transmission links, one that spreads an infection with high probability and one that does so only with a very small probability. It would then be natural to design a strategy that targets vertices with highly infectious connections.

Assume that $W\in\{a,b\}$ where $a<b$ and write $p_a=\PP(W=a)$ and $p_b=\PP(W=b)$. The strategy is defined so that each vertex is sampled independently with probability $s$ and, for each sampled vertex $u$, its neighbors with weight $b$ on their edge from $u$ are vaccinated: Recall that $\mathcal{N}_u$ denotes the set of neighbors of a vertex $u$ and let
$$
\mathcal{N}_u^{(b)}=\{v\in\mathcal{N}_u:W_{(u,v)}=b\},
$$
that is, $\mathcal{N}_u^{(b)}$ is the set of neighbors of $u$ for which the weight on the edge $(u,v)$ attains the larger value $b$. Then, if $u$ is sampled, the vertices in $\mathcal{N}_u^{(b)}$ are vaccinated. No action is taken if $\mathcal{N}_u^{(b)}$ is empty.

To derive the vaccination coverage, note that the probability that a randomly chosen vertex in the graph is not chosen for vaccination by a given neighbor equals
$$
\alpha=1-sp_b.
$$
As in the previous section we obtain the vaccination coverage from the equation
$$
1-v(s)=\sum_{j\geq 0}\alpha^jp_j.
$$
The derivation of the epidemic threshold is based on the same branching process as in the previous section, that is, an individual in the branching process consists of an unvaccinated vertex $u$ along with an outgoing edge $(u,w)$ that is not used for vaccination and that is open for transmission. To find an expression for the reproduction mean $R_s^{\verb"D"}$, which serves as the epidemic threshold, first note that in this case the degree distribution of vertex $w$ is not affected by the information that $w$ did not chose $u$ for vaccination (recall that the latter event is denoted $A$). Indeed, whether $u$ is vaccinated or not if $w$ is sampled is determined only by $W_{(w,u)}$. Hence $\PP_A(D_w=k)=\tilde{p}_k$.

Conditionally on $D_w=k$, the probability that $w$ is not vaccinated via any of its other $k-1$ neighbors (apart from $u$) is given by $\alpha^{k-1}$. We also need to determine the expected number of dangerous edges from $w$ to vertices in $\mathcal{N}_w\setminus \{u\}$ conditionally on that $D_w=k$ and on $A$ (note that, conditionally on the degree, the distribution of the number of dangerous edges from $w$ is not affected by the information that $w$ is not vaccinated). For this we need the corresponding probability that $w$ is sampled to name a neighbor for vaccination. With $V_w\in\{0,1\}$ denoting the number of times that $w$ is sampled to name a neighbor, we get
$$
\PP_{A,k}(V_w=1)=\frac{\PP_k(A|V_w=1)\PP_k(V_w=1)}{\PP_k(A)}=\frac{p_as}{\alpha}=:\nu.
$$
Note that this probability does not depend on $k$. If $V_w=0$, then the expected number of dangerous edges from $w$ (to other vertices than $u$) is $(k-1)\gamma$. If $V_w=1$ on the other hand, then the neighbors reached by edges with the large weight are vaccinated. The expected number of remaining out-edges from $w$ (to other vertices than $u$) is $(k-1)p_a$ and each one of these is open for transmission with probability $a$. The expected number of  dangerous edges from $w$ is hence $(k-1)ap_a$. Write $R_s^{\verb"W2"}$ for the basic reproduction number with the current vaccination strategy. We get
\begin{equation}\label{eq:RsW2}
R_s^{\verb"W2"}=(\nu ap_a+(1-\nu)\gamma)\sum_{k\geq 2}\tilde{p}_k\alpha^{k-1}(k-1).
\end{equation}
We now compare this to the epidemic threshold (\ref{eq:RvU}) for uniform vaccination and, in particular, to the threshold (\ref{eq:RbA}) for the standard acquaintance vaccination strategy.\medskip

\noindent\textbf{Example 4.2.1.} Figure 6 shows the basic reproduction numbers when the degree distribution is Po(14) and the weight distribution is specified by $\PP(W=0.1)=1-\PP(W=1)=0.9$ (most edges hence have a very small weight, but a small fraction has weight 1, implying almost sure transmission). The plot reveals that the weight based strategy clearly outperforms the other strategies in this case. The critical vaccination coverage is lowered from 0.58 with standard acquaintance vaccination to 0.48 with the weight based strategy.\hfill$\Box$\medskip

\noindent\textbf{Example 4.2.2.} Figure 7 shows the basic reproduction numbers for the same weight distribution as in the previous example when the degree distribution is a power-law with exponent 3.5 and mean 14. Again we see that the weight based strategy is the most efficient. \hfill$\Box$\medskip

\noindent\textbf{Example 4.2.3.} Finally, Figure 8 displays the basic reproduction numbers for the same power-law degree distribution as in the previous example but for a weight distribution specified by $\PP(W=0.1)=1-\PP(W=1)=0.5$. In this case almost nothing is gained by using the weight based strategy compared to standard acquaintance vaccination (the lines are almost aligned). The explanation for this is that, although the weight based strategy targets highly infective links, it does so more ``locally'' in the graph: Recall that \emph{all} neighbors with large weight on their edges from a sampled vertex are vaccinated. This means that, to achieve a given vaccination coverage, a much smaller sample of vertices is required compared to standard acquaintance vaccination if the probability of the larger weight is reasonably large; Figure 9 shows a plot for the current example. Thus the weight based strategy affects fewer parts of the graph and this cancels the positive effect that lies in securing high risk connections. However, the strategy does not perform worse than the standard acquaintance strategy. Hence the strategy is still more effective in the sense that it requires a smaller sample of vertices to name neighbors for vaccination to obtain a given vaccination coverage. In situations when there are costs associated with selecting and communication with the sampled vertices, this might be important.\hfill$\Box$\medskip

\begin{figure}[p]
\enlargethispage*{1cm}
\begin{center}
\mbox{\epsfig{file=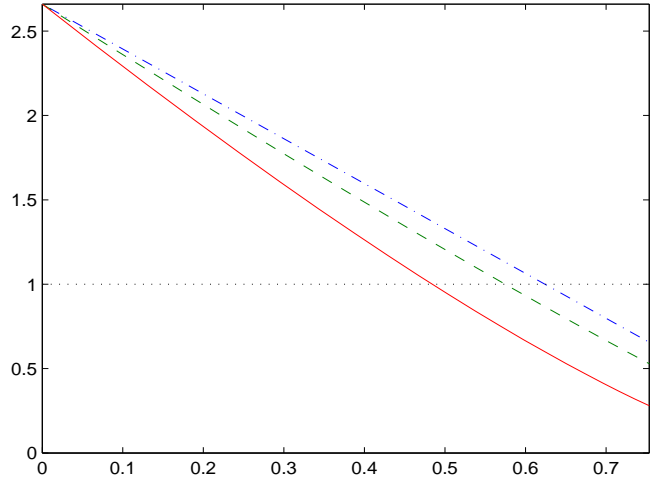,width=0.8\textwidth,
height=0.35\textheight}}
\end{center}
\caption{Basic reproduction numbers plotted against the vaccination coverage with $\PP(W=0.1)=1-\PP(W=1)=0.9$ and a Po(14) degree distribution: the weight based acquaintance strategy (solid line), the standard acquaintance vaccination (dashed line) and uniform vaccination (dash-dotted line).}
\end{figure}

\begin{figure}[p]
\begin{center}
\mbox{\epsfig{file=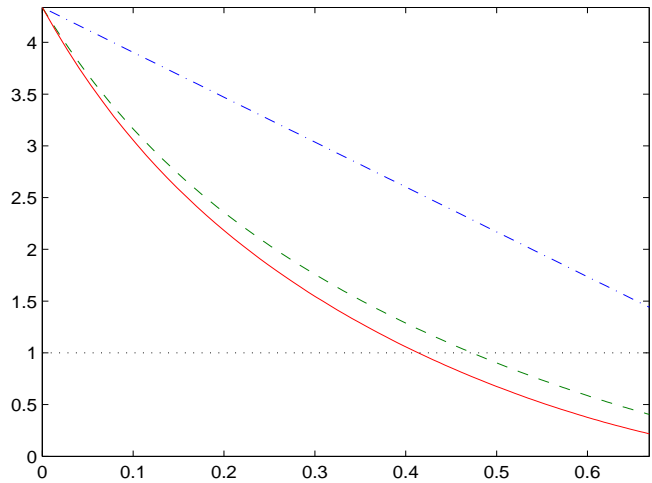,width=0.8\textwidth,
height=0.35\textheight}}
\end{center}
\caption{Basic reproduction numbers plotted against the vaccination coverage with $\PP(W=0.1)=1-\PP(W=1)=0.9$ and a power law degree distribution with exponent 3.5 and mean 14.}
\end{figure}

\begin{figure}[p]
\begin{center}
\mbox{\epsfig{file=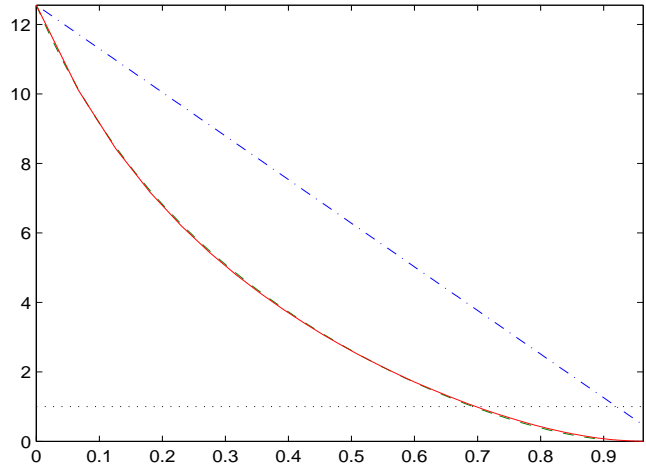,width=0.8\textwidth,
height=0.34\textheight}}
\end{center}
\caption{Basic reproduction numbers plotted against the vaccination coverage with $\PP(W=0.1)=1-\PP(W=1)=0.5$ and a power law degree distribution with exponent 3.5 and mean 14: the weight based acquaintance strategy (solid line), the standard acquaintance vaccination (dashed line) and uniform vaccination (dash-dotted line).}
\end{figure}

\begin{figure}[p]
\begin{center}
\mbox{\epsfig{file=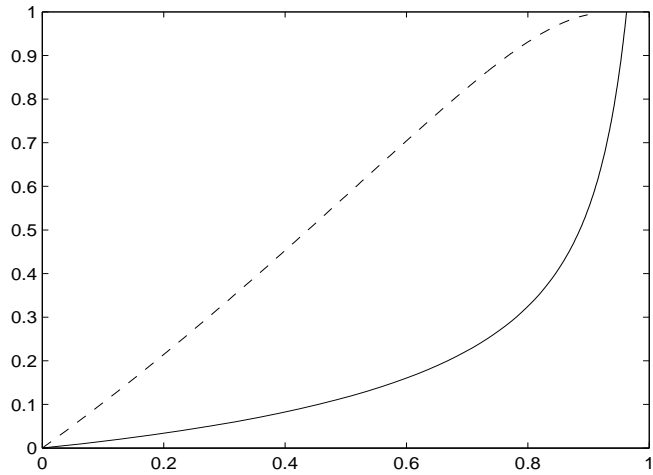,width=0.8\textwidth,
height=0.34\textheight}}
\end{center}
\caption{The fraction of the population that has to be sampled to vaccinate (at least one) neighbor(s) plotted against the resulting vaccination coverage for the weight based strategy (solid line) and standard acquaintance vaccination (dashed line).}
\end{figure}

\section{Summary and discussion}

We have formulated and analyzed a model for epidemic spread on weighted graphs, where the weight of an edge indicates the probability that it is used for transmission. Expressions have been derived for the epidemic threshold, specifying when there is a positive probability for an epidemic to take off. The case with independent weights is analogous to the case with a constant infection probability given by the mean weight. For degree dependent out-weights -- which for instance makes it possible to model a situation where high degree vertices infect their neighbors with a smaller probability -- however the behavior is different from a homogeneous epidemic.

Furthermore, we have analyzed a version of the acquaintance vaccination strategy where neighbors of the sampled vertices reached by edges with large weights are vaccinated. The selected vertices hence impose vaccination on the neighbor(s) that they have the strongest connection(s) to instead of a random neighbor. Two versions of this strategy have been treated: one for continuous weight distributions and one for two-point distributions. In the examples we have looked at, these strategies have been seen to outperform standard acquaintance vaccination, the difference being largest in cases where the weight distribution is highly right-skewed. The reason why the weight based acquaintance strategies perform better than standard acquaintance vaccination is that, in addition to removing the vaccinated neighbors, the ability to spread the epidemic is decreased also for the sampled vertices in that their high-weight connections are secured.

As for further work, there are numerous possibilities. In many situations it would be desirable to allow for (typically positive) correlations between the weights $W_{(u,v)}$ and $W_{(v,u)}$ on a given edge, for instance one might want to assign only one weight per edge, specifying the probability of transmission in any direction. This leads to complications in the current analysis, basically because the information that a vertex is unvaccinated then gives information on the weights on the edges of its neighbors. Furthermore, the basic idea in acquaintance vaccination is that, by vaccinating neighbors of the sampled vertices, one reaches vertices with higher degree. A natural further development of this idea would be to vaccinate neighbors with maximal degree, that is, selected vertices are asked to identify their neighbor(s) with the largest degree among the neighbors (assuming that they have this information) and these neighbors are then vaccinated. Unfortunately this seems to lead to complicated dependencies in the resulting epidemic process.

We also mention that it would be interesting to investigate the final size of the epidemic. This is usually related to the probability of a large outbreak and quantified via an equation involving the generating function of the reproduction distribution. For the vaccination strategies that we have considered here, this equation would involve the distribution (\ref{eq:ord_stat}) of order statistics and is hence presumably complicated. But it would be interesting to study the final size by aid of simulation. Other possible continuations include investigating how the results are affected by introducing clustering (triangles and other short cycles) in the underlying graph, to involve time-dynamic in the vaccination procedure and to generalize the model for the epidemic spread.\bigskip

\noindent\textbf{Acknowledgement.} The author gratefully acknowledges the support from The Bank of Sweden Tercentenary Foundation.

\end{document}